\documentstyle[a4,12pt]{article}

\newfam\msbfam
\font\twelmsb=msbm10 at 12pt
\font\tenmsb=msbm10
\font\sevenmsb=msbm10 at 7pt
\font\fivemsb=msbm10 at 5pt

\textfont\msbfam=\twelmsb
\scriptfont\msbfam=\sevenmsb
\scriptscriptfont\msbfam=\fivemsb

\let\msy\Bbb

\begin{document}
% latex, version january 1995
\newtheorem{definition}[subsubsection]{Definition}
\newtheorem{theorem}[subsubsection]{Theorem}
\newtheorem{lemma}[subsubsection]{Lemma}
\newtheorem{cortje}[subsection]{Corollary}
\newtheorem{lemmatje}[subsection]{Lemma}
\newtheorem{proposition}[subsubsection]{Proposition}
\newtheorem{corollary}[subsubsection]{Corollary}
\def\square{\Box}
%\newremark{example}[definition]{Example}
%\newremark{examples}[definition]{Examples}
\newtheorem{remark}[subsubsection]{Remark}
\newtheorem{observation}[subsubsection]{Observation}
%\newremark{assumption}[definition]{Assumption}
\newenvironment{prf}[1]{\trivlist
\item[\hskip
\labelsep{\it #1.\hspace*{.3em}}]}{~\hspace{\fill}~$\square$\endtrivlist}
\newenvironment{proof}{\begin{prf}{Proof}}{\end{prf}}

%MACROS----------------
\def\et{{\acute et}}
\def\Et{{\rm \acute Et}}
\def\Homext{{\cal H}om^{ext}(\pi_1(X/S),G)}
\def\L{\msy L} \def\cC{{\cal C}}
\def\cGMg{{}_G{\cal M}_g}             \def\GMg{{}_GM_g}
\def\bcGMg{\overline{{}_G{\cal M}_g}} \def\bGMg{\overline{{}_GM_g}}
\def\cMg{{\cal M}_g}                  \def\Mg{M_g}
\def\bcMg{\overline{{\cal M}_g}}       \def\bMg{\overline{M_g}}
\def\cMgan{{\cal M}^{an}_g}           \def\Mgan{M_g^{an}}
\def\bcMgan{\overline{{\cal M}_g}^{an}}   \def\bMgan{\overline{M_g}^{an}}
\def\cGMgan{{}_G{\cal M}^{an}_g}      \def\GMgan{{}_GM_g^{an}}
\def\bcGMgan{\overline{{}_G{\cal
M}_g}^{an}}\def\bGMgan{\overline{{}_GM_g}^{an}}
\def\p{{\bar p}}
\def\pirel{\pi_1(X/S,s)}
\def\pifi{\pi_1(X_p,s(p))}
\def\pif#1{\pi_1(X_{#1}, s(#1))}
\def\pip#1#2{\pi_1^{\L}(#1,#2)}
\def\pipf{\pi_1^{\L}(X_\p, s(\p))}
\def\piprel{\pi_1^{\L}(X/S,s)}
\def\q{{\bar q}}
\def\mapright#1{\smash{\mathop{\longrightarrow}\limits^{#1}}}
\def\mapleft#1{\smash{\mathop{\longleftarrow}\limits^{#1}}}
\def\mapdown#1{\Big\downarrow\rlap{$\vcenter{\hbox{$\scriptstyle#1$}}$}}
\def\downmap#1{\downarrow\rlap{$\vcenter{\hbox{$\scriptstyle#1$}}$}}
\def\mapup#1{\Big\uparrow\rlap{$\vcenter{\hbox{$\scriptstyle#1$}}$}}
\newcommand\Gk{\mbox {$G^{(k)}$}}
\newcommand\Gkn {\mbox {$G^{(k),n}$}}
\newcommand\Gkeen {\mbox {$G^{(k+1)}$}}
\newcommand\pk {\mbox {$\pi^{(k)}$}}
\newcommand\pkn {\mbox {$\pi^{(k),n}$}}
\newcommand\pkeenn {\mbox {$\pi^{(k+1),n}$}}
\newcommand\pkeen {\mbox {$\pi^{(k+1)}$}}
\newcommand\Okn {\mbox {$O^{(k),n}$}}
\newcommand\Akn {\mbox {$A^{(k),n}$}}
\newcommand\Ikn {\mbox {$I^{(k),n}$}}
\newcommand\Ok {\mbox {$O^{(k)}$}}
\newcommand\Ak {\mbox {$A^{(k)}$}}
\newcommand\Ik {\mbox {$I^{(k)}$}}
\newcommand\Akeen {\mbox {$A^{(k+1)}$}}
\newcommand\Akeenn {\mbox {$A^{(k+1),n}$}}
\newcommand\ra{\rightarrow}

\title{Moduli of curves with non-abelian level structure}
\author{M. Pikaart \and A.J. de Jong\thanks{{\rm The research of
A.J.~de Jong has been made possible by a fellowship
of the Royal Netherlands Academy of Arts and Sciences.}}}

\maketitle

\section{Introduction}

Deligne and Mumford introduced the moduli stack
$\cGMg$ parametrizing smooth genus $g$ curves with Teichm{\"u}ller
structure of level $G$, a finite group. For example,
if $G=\{e\}$, resp.\ $G\cong (\msy Z/n\msy Z)^{2g}$, this moduli
stack is just the moduli stack of smooth curves $\cMg$, resp.\
of smooth curves with abelian level structure (sometimes
denoted $\cMg[n]$). They also defined $\bcMg$, the moduli
stack of stable curves and proved it is proper over
${\rm Spec}(\msy Z)$. Taking the normalization of $\bcMg\times
{\rm Spec}(\msy Z[1/\# G])$ in the function field
of $\cGMg$ defines $\bcGMg$, proper over ${\rm Spec}(\msy Z[1/\# G])$.

Let $\Pi$ denote the standard fundamental group of a compact
Riemann surface of genus $g$. The nth powers together with
the kth order commutators generate a normal subgroup $\Pi^{(k),n}$.
(We regard $[a,b]$ as a commutator of order 2.) Let $G$ be the
quotient $\Pi/\Pi^{(k+1),n}$. We show that if $k\geq 1$, $n\geq 3$
the coarse moduli scheme $\bGMg$ for $\bcGMg$ exists, and we actually
have $\bGMg\cong\bcGMg$. The following theorem is our main result.

\medskip
\noindent{\bf Theorem\enspace \ref{glad}} {\sl
Suppose $k\in \{1,2,3\}$ and $n\geq 3$. The structural morphism
$\bGMg\to {\rm Spec}(\msy Z[1/n])$ is smooth if and only if
\begin{itemize}
\item{$k=1$} and $g=2$,
\item{$k=2$} and $n$ is odd,
\item{$k=3$} and $n$ is odd or $n$ is divisible by $4$.
\end{itemize}
Furthermore, if $k \geq 4$, $n \geq 3$ and $n$ relatively prime to $6$,
then $\bGMg\to {\rm Spec}(\msy Z[1/n])$ is smooth.}
\medskip

Looijenga constructed a smooth and compact cover
$\overline{M_g[_2^n]}^{an},\ n\geq 3$ of $\bMgan$ using Prym level
structures. He proves that these coverings are universally ramified
along the boundary of $\cMgan$, see \cite{Looijenga}.
This holds also for our coverings
as can be seen from \ref{mono}. These coverings can be applied
to the construction of the Chow rings of the stacks $\bcMg[n]$:
for example it is clear that the specialization maps are ring homomorphisms.

The paper is organized as follows. Section 2 deals with the definition
of the moduli problem. The formulation is in terms of the relative
fundamental group of a family of curves.
Its final proposition states that $\bGMg\to {\rm Spec}(\msy Z[1/n])$
is smooth if and only if the associated compact analytic space $\bGMgan$ is
is a manifold. Section 3 contains a precise description of the monodromy
along the boundary of $\bGMgan$ for general groups $G$ and the statement
of the main result. Sections 4 and 5 are of a topological nature.
We describe the monodromy on the relative fundamental group
of the universal curve along the boundary in terms of Dehn twists.
We have to describe the situation in some detail in order to
understand how this monodromy acts on the finite quotients $\Pi/\Pi^{(k),n}$.
These considerations prove the main theorem, using Section 6,
which contains the necessary computations for a free group on three
generators.

The authors were stimulated by the article \cite{Looijenga}.
We thank Prof.~Looijenga for numerous discussions explaining
his and other results. We thank Prof.~Oort, who remarked that
it should be possible to do everything algebraically and
drew our attention to the article \cite{Oda}.

\subsection{Notations and conventions}\label{notation}
\begin{enumerate}
\item Throughout the paper $g$ is a fixed natural number at least 2.
\item $G$ is a finite group.
\item $\msy L$ is a set of primes, $\msy L$ contains the primes dividing
$\# G$.
\item The stack of stable curves of genus $g$ is denoted $\bcMg$, the
open substack of smooth curves $\cMg$. Stacks are denoted
by script letters. For definitions and results concerning stacks
we refer to \cite{DM}.
\item Suppose $\Gamma$ is a (pro-finite) group. A characteristic
subgroup of $\Gamma$ is a normal subgroup fixed by any automorphism
of $\Gamma$. A characteristic quotient is one whose kernel is a
characteristic subgroup. If $\Gamma$ is profinite and topologically
finitely generated then it is the direct limit of its finite
characteristic quotients.
\item Let $\Gamma$ denote a (finitely generated) group. For $a,b\in \Gamma$ we
put $[a,b]=a^{-1}b^{-1}ab$, so that $ab=ba[a,b]$. We define the lower
central series $\Gamma^{(k)}$ of $\Gamma$ by $\Gamma^{(1)}=\Gamma$ and
$\Gamma^{(k+1)}=[\Gamma^{(k)},\Gamma]$. The subgroup of $\Gamma$ generated
by nth-powers is denoted $\Gamma^n$. We write $\Gamma^{(k),n}$ to indicate
the subgroup generated by $\Gamma^{(k)}$ and nth-powers,
$\Gamma^{(k),n}=\Gamma^{(k)}\cdot \Gamma^n$. Any group homomorphism
maps commutators to commutators and nth-powers to nth-powers, hence
preserves these subgroups. In particular, the subgroups $\Gamma^{(k),n}$
are characteristic subgroups of $\Gamma$.
\item Let $\Pi=\Pi_g$ denote the standard fundamental group $\Pi=\pi_1(S)$
of a compact Riemann surface $S$ of genus $g$. \end{enumerate}

\section{Definition of the moduli problem}

In this section we recall the definition of the moduli problem of
Teichm\"uller level structures, see \cite[ Section 5]{DM}. Furthermore,
we prove that the Deligne-Mumford compactification of the
associated stack is smooth if and only if the corresponding
analytic orbifold is smooth.

\subsection{The relative fundamental group}

In this section we define the relative fundamental
group for a proper smooth morphism $f: X\to S$ with connected
fibres and endowed with a section $s: S\to X$. In order to
motivate the definition in the algebraic case (and since we
need it also) we first do the analytic case.

\subsubsection{The analytic case}

Here $f:X\to S$ is a proper smooth morphism of analytic spaces
with connected fibres. In addition we are given a section $s:S\to X$
of $f$. We define a locally constant sheaf of groups
$\pirel$ over $S$ such that for all points $p\in S$ we have an isomorphism
of groups
$$\pirel_p\cong \pifi,$$
of the fibre of the sheaf at $p$ with the topological fundamental group
of the fibre $X_p$ of $f$ at $p$ with base point $s(p)$.

To construct $\pirel$ we choose for any point $p$ of $S$ a connected open
neighbourhood $U_p\subset S$ and a topological isomorphism
$$\phi_p: f^{-1}(U_p)\cong X_p\times U_p.$$
Such can be found compatible with $f$ and the projection to $U_p$,
inducing the identity on $X_p$ and
such that $\phi_p\circ s$ equals $q\mapsto (s(p), q)$. Over $U_p$ we
take $\pirel$ constant with fibre $\pifi$.

To glue these we note that given two points $p_1, p_2\in S$ there is for
any $q\in U_{p_1}\cap U_{p_2}$ an identification
$$X_{p_1}=X_{p_1}\times \{q\}
\mapright{\phi_{p_1}^{-1}} X_q \mapright{\phi_{p_2}}
X_{p_2}\times \{q\}=X_{p_2}.$$
This identification is compatible with base points $s(p_i)$ and
depends continuously on $q\in U_{p_1}\cap U_{p_2}$. This means that
the induced isomorphism
$$ \pif{p_1}\cong \pif{p_2}$$
is constant on the connected components of $U_{p_1}\cap U_{p_2}$. Hence we
get the desired gluing. We leave to the reader the trivial verification
that these gluings satisfy the desired cocycle condition on
$U_{p_1}\cap U_{p_2}\cap U_{p_3}$.

If we choose other $\phi_p$, say $\phi'_p$, then for $q\in U_p$ the map
$X_q\to X_p\times \{q\}\to X_q$, using first $\phi_p^{-1}$ then $\phi'_p$,
is homotopic to the identity. Hence the resulting sheaves $\pirel$
are canonically isomorphic. A similar argument deals with the shrinking
of the neighbourhoods $U_p$.

\begin{proposition}The construction given above defines a locally constant
sheaf of groups $\pirel$ over $S$, characterized by the following
properties:
\begin{enumerate}
\item For any point $p\in S$
there is given an isomorphism $\pirel_p\cong \pifi$.
\item The monodromy action of $\pi_1(S,p)$ on the fibre of the locally
constant sheaf
$$\rho : \pi_1(S,p)\to {\rm Aut}\big(\pirel_p\big) $$
agrees, via the isomorphism of 1), with the action of $\pi_1(S,p)$
on $\pifi$ deduced from the split exact sequence
$$1\longrightarrow \pifi \longrightarrow \pi_1(X, s(p))
\lower3pt\hbox{$\longrightarrow$}\llap{\raise3pt\hbox{$\mapleft{s_\ast}$}}
\pi_1(S,p)\longrightarrow 1.$$
\item The construction of $\pirel$ commutes with arbitrary base change
$S'\to S$.
\end{enumerate}\end{proposition}

\begin{proof}Suppose $\gamma\in \pi_1(S,p)$ and $\alpha\in \pifi$. It is
well known and easy to prove that $s_\ast(\gamma)\alpha s_\ast(\gamma^{-1})$
is equal to the horizontal transportation of the loop $\alpha$ over
$\gamma$. Clearly this describes the monodromy representation for the sheaf
$\pirel$. The proof of the other assertions is left to the reader.\end{proof}

\begin{remark}\label{innerconjugation}
{\rm Suppose $s':S\to X$ is a second section of $f$. In
general $\pirel$ is not isomorphic to $\pi_1(X/S,s')$. However, locally
on $S$, say over $U\subset S$, we can choose a homotopy $H$ between
$s$ and $s'$. This will induce an identification
$$ i_H: \pirel|_U\cong \pi_1(X/S,s')|_U.$$
This is unique up to an {\it inner automorphism} of $\pirel$. Indeed,
if $H'$ is another such homotopy, then combining $H$ and $H'$ gives a
familly of loops in $X$ over $U$, with base points $s(u)$, i.e., a section
of $\pirel$ over $U$. The map $(i_{H'})^{-1}\circ i_H$ is equal to
conjugation with this section.}\end{remark}

\subsubsection{The algebraic case}

Here we consider a proper smooth morphism of schemes
$f:X\to S$ with connected geometric fibres. As before we have a section
$s:S\to X$ of the morphism $f$. Further, we assume given a set of primes
$\L$ such that all residue characteristics of $S$ are {\it not} in $\L$.

We recall some general notations concerning algebraic fundamental groups.
We refer to \cite{SGA1} and \cite{Murre} for more details.
If $Y$ is a scheme, then $\Et(Y)$ denotes the category of finite \'etale
coverings $Y'\to Y$. If $\p$ is a geometric point of $Y$ then we denote
by $F_{Y,\p}$ the fundamental functor (or fibre functor)
$ F_{Y,\p} : \Et(Y)\longrightarrow Set $
which associates to $Y'$ the set of geometric points $\p\to Y'$ lying over
$\p\to Y$. By definition we have $\pi_1(Y,\p)={\rm Aut}(F_{Y,\p})$; this is
the fundamental group of $Y$ with base point $\p$. To compare the fundamental
groups with base points $\p$, resp.\ $\q$ we use a path from $\p$ to $\q$,
i.e., an isomorphism of fibre functors
$ \alpha: F_{Y,\p}\longrightarrow F_{Y,\q}.$
Obviously, $\alpha$ gives an isomorphism
$\alpha_\ast : \pi_1(Y,\p)\longrightarrow \pi_1(Y,\q).$
We note that it is independent of the choice of $\alpha$ up to conjugation.
A morphism of schemes $h: Y\to Z$, defines a functor
$h^\ast : \Et(Z)\to \Et(Y)$, which satisfies $F_{Y,\p}=F_{Z,h(\p)}\circ
h^\ast$.
Therefore we get $h_\ast$ on loops and on paths.

A slight modification of the above gives $\pip Y\p$, the algebraic fundamental
group classifying Galois coverings of degree in $\L$. Formally it can be
defined as
$$\pip Y\p = \lim_{{\longleftarrow}} G, $$
where the limit is taken over all surjections $\pi_1(Y,\p)\to G$ onto finite
groups $G$ whose orders have only prime factors from $\L$.

In the sequel we will use the following results from \cite{SGA1}: The
sequence
$$\pi_1(X_\p,s(\p))\longrightarrow \pi_1(X,s(\p)) \longrightarrow
\pi_1(S,\p)\longrightarrow 1 $$
is exact. If we take the pushout of this sequence with the surjection
$\pif\p\to\pip {X_\p}\p$ then the sequence also becomes left exact
$$ 1\longrightarrow \pip {X_\p}\p \longrightarrow\pi'_1(X,s(\p))
\lower3pt\hbox{$\longrightarrow$}\llap{\raise3pt\hbox{$\mapleft{s_\ast}$}}
\pi_1(S,\p)\longrightarrow 1.\eqno{(*)}$$
See \cite[Expos\'e XII 4.3, 4.4]{SGA1}.
As before the section $s$ defines a splitting
$s_\ast : \pi_1(S,\p)\to \pi_1(X,s(\p))$.

\begin{proposition}(\cite[Expos\'e XII 4.5]{SGA1})\label{piprel}
There is a pro-object in the category of locally
constant sheaves of groups on $S_\et$, denoted $\piprel$,
determined up to unique isomorphism by the following properties:
\begin{enumerate}
\item For any geometric point $\p$ of $S$ there is given an
isomorphism $$\piprel_\p\longrightarrow \pip {X_\p}{s(\p)}.$$
\item The monodromy presentation
$$\rho : \pi_1(S,\p)\longrightarrow {\rm Aut}\big(\piprel_\p\big)$$
equals, via 1), the action of $\pi_1(S,\p)$ on $\pip {X_\p}{s(\p)}$
deduced from {\rm (*)}.
\item The construction of $\piprel$ commutes with arbitrary base change
$S'\to S$.\end{enumerate}\end{proposition}

\begin{proof}Of course we may assume that $S$ is connected. Take a
geometric point $\p$ of $S$. First we note that, since $\pipf$ is
topologically finitely generated, it is the direct limit of its
characteristic finite quotients:
$$\pipf=\lim_{{\scriptstyle \longleftarrow}\atop {\scriptstyle \omega}}
G_\omega$$
The action of $\pi_1(S,\p)$ on $\pipf$ deduced from (*)
gives an action $\rho_\omega$ on each $G_\omega$. This defines a finite locally
constant \'etale sheaf ${\cal F}_\omega$ on $S_\et$ whose fibre in $\p$
is given by $G_\omega$ and monodromy action equal to $\rho_\omega$.
We put
$$ \piprel=\lim_{{\scriptstyle \longleftarrow}\atop {\scriptstyle \omega}}
{\cal F}_\omega.$$
This immediately gives 1) and 2) for our chosen point $\p$. Part 3) is also
clear if there exists a lift of $\p$ to a geometric point of $S'$, see
\cite{SGA1}.

Thus it suffices to prove 1) and 2) for a second geometric point $\q$ of $S$.
Note that by definition
$$\pi_1(X,s(\p))={\rm Aut}\big(F_{S,\p}\circ s^\ast\big).$$
Hence if we choose an isomorphism $\alpha : F_{S,\p}\to F_{S,\q}$
then we get a commutative diagram
$$\matrix{\pi_1(X, s(\p))&
\lower3pt\hbox{$\longrightarrow$}\llap{\raise3pt\hbox{$\mapleft{s_\ast}$}}&
\pi_1(S,\p)\cr
\mapdown{\alpha_\ast}&&\mapdown{\alpha_\ast}\cr
\pi_1(X, s(\q))&
\lower3pt\hbox{$\longrightarrow$}\llap{\raise3pt\hbox{$\mapleft{s_\ast}$}}&
\pi_1(S,\q).\cr}$$
By (*) this induces an isomorphism
$$ \alpha^{\L} : \pipf \longrightarrow \pi_1^{\L}(X_\q, s(\q)).$$
This is an isomorphism compatible with the actions of
$\pi_1(S,\p)$ and $\pi_1(S,\q)$, compared via $\alpha_\ast$.

The fibre of $\piprel$ at $\q$ is by definition
$$F_{S,\q}\big(\piprel\big)\mapright{\alpha^{-1}}
F_{S,\p}\big(\piprel\big)= \pipf.$$
If we use $\alpha^{\L}$ to identify this with $\pi_1^{\L}(X_\q, s(\q))$ then
we see by the above that the monodromy action on this exactly corresponds
to the action deduced from (*) (with $\q$ in stead of $\p$).
We leave to the reader the verification that another choice of $\alpha$
gives the same identification
$$\piprel_\q\longrightarrow \pi_1^{\L}(X_\q, s(\q)).$$\end{proof}

Suppose $s': S\to X$ is a second section of $f$. Take a geometric point
$\p$ of $S$. We write $i_\p$ for the morphism $X_\p\to X$.
We say that a path $\beta$ on $X$ connecting $s(\p)$ to
$s'(\p)$ lies in $X_\p$ if there exists a path $\tilde \beta$ in $X_\p$
such that $i_{\p,\ast}(\tilde \beta)=\beta$. We remark that any path
$\beta$ connecting $s(\p)$ to $s'(\p)$ lies in $X_\p$ if and only if
$f_\ast(\beta)=1$ (in $\pi_1(S,\p)$). This is easily seen using the
first exact sequence above.

Let us take such a $\beta$ lying in $X_\p$. It gives rise to a commutative
diagram
$$\matrix{1&\longrightarrow&
\pipf&\longrightarrow&\pi'_1(X,s(\p))&\longrightarrow&\pi_1(S,\p)&
\longrightarrow&1\cr
&&\mapdown{\cong}&&\mapdown{\beta_\ast}&&\mapdown{{\rm id}}&&\cr
1&\longrightarrow&\pi_1^{\L}(X_\p, s'(\p))
&\longrightarrow&\pi'_1(X,s'(\p))&\longrightarrow&\pi_1(S,\p)&
\longrightarrow&1\cr}$$
The isomorphism $\pipf\cong \pi_1^{\L}(X_\p, s'(\p))$ determined in this way
is unique up to inner conjugation. We constructed $\piprel$, resp.\
${\pi_1^{\L}(X/S,s')}$ as the limit of sheaves ${\cal F}_\omega$, resp.\
${\cal F}'_\omega$ corresponding to characteristic quotients
$\pipf\to G_\omega$, resp.\ ${\pi_1^{\L}(X/S,s')}\to G_\omega$. Note
that the isomorphisms $G_\omega\to G'_\omega$ induced from the above
are also unique up to inner conjugation.

Let $S_\omega\to S$ be the finite \'etale covering of $S$ trivializing
the action of $\pi_1(S,\p)$ on both $G_\omega$ and $G'_\omega$. We can
use the above to get an isomorphism of (constant) sheaves of groups
$$ {\cal F}_\omega|_{S_\omega}\longrightarrow {\cal F}'_\omega|_{S_\omega}.$$
We claim this isomorphism is unique up to inner conjugation.
This is clear if we only change $\beta$, but what happens if we change $\p$
to $\q$? Take $\alpha : F_{S,\p}\to F_{S,\q}$ as in the proof
of Proposition \ref{piprel}. What we have to check is that
$s'_\ast(\alpha)\circ\beta\circ s_\ast(\alpha^{-1})$ is a path
connecting $s(\q)$ to $s'(\q)$ lying in $X_\q$. But this is clear
since $f_\ast\big(s'_\ast(\alpha)\circ\beta\circ s_\ast(\alpha^{-1})\big)
=\alpha\circ \alpha^{-1}=1.$

\begin{corollary}\label{invariance}The construction above defines
locally in the \'etale topology on $S$ identifications of the
finite quotients of the sheaves $\piprel$ and ${\pi_1^{\L}(X/S,s')}$.
These identifications are unique up to inner conjugation
and agree via 1) of Proposition \ref{piprel} with the usual
identifications of $\pip {X_\p}{s(\p)}$ and
$\pip {X_\p}{s'(\p)}$.\end{corollary}

\subsection{Exterior homomorphisms}

In this section we define the sheaf of exterior homomorphisms
of the relative fundamental group of $X$ over $S$ into a fixed finite
group $G$. As in \cite[5.5]{DM} this sheaf will be denoted $\Homext$ and will
be a finite locally constant sheaf of sets on $S$ (or $S_\et$).
We note that $\pi_1(X/S)$ has not been defined.

\subsubsection{The analytic case}

Here $f:X\to S$ is a proper smooth morphism of analytic spaces
with connected fibres.

If $f$ has a section $s$ then we can look at the locally constant
sheaf
$$ {\cal F}={\cal H}om(\pirel, G)$$
on $S$. It has finite fibres since $\pifi$ is finitely generated for all
$p$ in $S$. There is a natural action of the sheaf $\pirel$ on the
sheaf ${\cal F}$ given by conjugation. We define the sheaf of exterior
homomorphisms as the quotient of ${\cal F}$ by this action
$$ \Homext:={\cal F}\big/\pirel .$$
It is clear from Remark \ref{innerconjugation} that the right hand side
does not depend on the chosen section $s$.

In general, we choose an open covering $S=\bigcup U_i$ such that
$X\to S$ has a section over each $U_i$. The sheaf $\Homext$ is then
defined by gluing the sheaves constructed above.
The fibres are described by the formula:
$$\Homext_p={\rm Hom}(\pi_1(X_p,q),G)\big/ \pi_1(X_p,q), $$
where $q$ is any point of the fibre $X_p$. Note that the monodromy
action of $\gamma\in \pi_1(S,p)$ on this is given by horizontal
transport of loops in $\pi_1(X_p,q)$.

\subsubsection{The algebraic case}

Here we assume that $S$ is a scheme over ${\rm Spec}(\msy Z[1/\#G])$.
The morphism $f:X\to S$ is still assumed proper smooth with
connected geometric fibres.

The construction of $\Homext$ in this case is exactly the same as
for the analytic case. We use $\piprel$, where $\L$ is the set of
primes dividing $\# G$, and we use that sections of $f$ exist
locally in the \'etale topology on $S$. We use also Corollary
\ref{invariance}.
The geometric fibres of the resulting sheaf can be described as follows:
$$\Homext_\p={\rm Hom}(\pi^{\L}_1(X_\p,\q),G)\big/ \pi^{\L}_1(X_\p,\q)=
{\rm Hom}(\pi_1(X_p,q),G)\big/ \pi_1(X_p,q), $$
the last equality holds in view of our definition of $\L$.
This justifies dropping $\L$ from the notations.

\subsubsection{Comparison}

We note that if $S$ is a scheme of finite type over $\msy C$
there is a canonical homomorphism
$\pi_1(X^{an}/S^{an},s)\to \piprel^{an}$, identifying the relevant
finite quotient sheaves. Clearly this gives rise to an identification
of sheaves of exterior homomorphisms.

\subsection{Teichm\"uller level structures}

In this section the morphism $f:X\to S$ will be a familly of smooth
projective curves of genus $g$. The abstract finite group $G$ will
be fixed. In both the analytic case and the algebraic case we make the
following definition.

\begin{definition}{\rm \cite[5.6]{DM}} A Teichm\"uller
structure $\alpha$ of level $G$ on $X\to S$
is a surjective exterior homomorphism
$$\alpha\in \Gamma(S, \Homext).$$
Thus locally on $S$ (resp.\ $S_\et$) $\alpha$ corresponds to a
surjective homomorphism $\pirel\to G$.\end{definition}

We want to consider the moduli spaces parametrizing
smooth stable curves of genus $g$ with a Teichm\"uller structure
of level $G$. However, as usual, it is more convenient to work
with stacks. Thus let $\cGMg$ denote the stack whose
category of sections over the scheme $S$
(lying over ${\rm Spec}(\msy Z[1/\#G])$) is the category of smooth stable
curves $X\to S$ of genus $g$ endowed with a Teichm\"uller structure
of level $G$, see \cite[Section 5]{DM}. The construction of $\Homext$
shows that the stack $\cGMg$ is representable
finite \'etale over $\cMg[1/\#G]$. Thus
$\cGMg$ is a separated algebraic stack, smooth over
${\rm Spec}(\msy Z[1/\#G])$.

In a similar way we define the analytic stack
$\cGMgan$ classifying complex curves with a Teichm\"uller structure of
level $G$. By comparing the stacks $\cGMg$ and $\cGMgan$ with the
stack $\cMg$ and using the known result for $M_g$ one derives easily
the following result.

\begin{theorem}\label{existence}With notations as above.
\begin{enumerate}
\item A coarse moduli scheme ${}_GM_g$ for $\cGMg$ exists. It is separated
of finite type over ${\rm Spec}(\msy Z[1/\#G])$.
\item A coarse analytic moduli space ${}_GM_g^{an}$ for $\cGMgan$
exists. It is canonically isomorphic to the analytic space associated
to the complex variety ${}_GM_g\otimes \msy C$.
\end{enumerate}\end{theorem}

Suppose we have a surjection $G\to G'$. There is a natural forgetful
morphism $\cGMg\to {}_{G'}\cMg$. This morphism is representable
finite \'etale. Hence, if ${}_{G'}\cMg$ is a scheme (i.e., isomorphic
to ${}_{G'}M_g$), then so is $\cGMg$. In this case the morphism
${}_{G'}M_g\to \GMg$ is finite \'etale (perhaps of degree $0$).
Finally, we note that $\cGMg$ might be empty, this being the case
if $G$ is not isomorphic to a quotient of the fundamental group of
a Riemann surface of genus $g$.

\subsubsection{Abelian level structures}

In this section we treat the algebraic case. Let $f:X\to S$
be a smooth stable curve of genus $g$ over $S$. Fix a natural number
$m$. An abelian structure structure of level $m$ on $X$ over $S$ is
defined as an isomorphism of \'etale sheaves over $S$
$$ (\msy Z/m\msy Z)^{2g}_S\longrightarrow R^1f_\ast(\msy Z/m\msy Z).$$
Such a level structure can only exist if the base scheme $S$ lies
over ${\rm Spec}(\msy Z[1/m])$; let us assume this is the case.

Note that the sheaf ${\cal H}om^{ext}(\pi_1(X/S),\msy Z/m\msy Z)$
is isomorphic to the subsheaf of primitive elements
in $R^1f_\ast(\msy Z/m\msy Z)$.
Thus we see that the moduli stack of curves with an abelian level $m$
structure is isomorphic to the stack $\cGMg$
with $G=(\msy Z/m\msy Z)^{2g}$.
In particular, if $m\geq 3$, then $\GMg$ is a fine moduli
scheme smooth over ${\rm Spec}(\msy Z[1/m])$. The following is deduced
from the above.

\begin{proposition}
If the finite group $G$ allows a surjection onto $(\msy Z/m\msy Z)^{2g}$
for some $m\geq 3$ then the coarse moduli scheme ${}_GM_g$ is a
fine moduli scheme.\end{proposition}

\subsubsection{Compactifications}

In order to get compact moduli spaces we just take the normalization
with respect to the Deligne-Mumford compactification. In this subsection
the finite group $G$ will be fixed, of order $n=\#G$.

Consider the Deligne-Mumford compactification $\cMg\subset \bcMg$.
We define $\bcGMg$ as the normalization of $\bMg[1/n]$ with
respect to $\cGMg$. Similarly we define $\bcGMgan$ as the
normalization of $\bcMgan$ with respect to $\cGMgan$.
As per convention we denote $\bGMg$ (resp.\ $\bGMgan$) the associated
coarse moduli scheme (resp.\ analytic moduli space).
The morphism
$$\bGMg\longrightarrow {\rm Spec}(\msy Z[1/n])$$
is proper, since $\bMg[1/n]$ is proper over ${\rm Spec}(\msy Z[1/n])$.
To see whether this morphism is smooth we have the following
criterion.

\begin{proposition}\label{criterium}
For finite groups $G$ which allow a surjection
$G\to (\msy Z/m\msy Z)^{2g}$ for some $m\geq 3$ the following
statements are equivalent:
\begin{enumerate}
\item The morphism $\bGMg\longrightarrow {\rm Spec}(\msy Z[1/n])$ is smooth.
\item The analytic space $\bGMgan$ is a (nonsingular) complex manifold.
\end{enumerate}\end{proposition}

\begin{proof}We use that $\bGMgan$ is isomorphic to the analytic space
associated to the variety $\bGMg\otimes \msy C$. Thus it suffices
to show that $\bGMg\otimes \msy C$ nonsingular implies
$\bGMg\otimes \bar {\msy F}_p$ nonsingular, where $p>0$ and $p$ does not
divide $n$. The argument will be based on the fact that the morphism
$\varphi : \bGMg\to \bMg[1/n]$ is tamely ramified along the boundary.

To see this we need a description of the complete local rings
of $\bGMg$ in points on the boundary. Suppose that $C$ is a stable curve
of genus $g$ over an algebraically closed field $k$ of characteristic $p$.
The singular points of $C$ are $P_1,\ldots, P_\ell$.
Let $\cC\to {\rm Spf}\big(W(k)[[t_1,\ldots,t_{3g-3}]]\big)$ be the universal
deformation of $\cC$. We choose the parameters $t_i$ such that $t_i=0$
defines the locus where $P_i$ survives as a singular point,
for $i=1,\ldots, \ell$.
Put $A=W(k)[[t_1,\ldots,t_{3g-3}]]$. Since $\cC$ is (uniquely) algebraizable
we have a morphism ${\rm Spec}(A)\to \bcMg$ and
$${\rm Spec}\big(A[1/t_1\ldots t_\ell]\big)\longrightarrow \cMg.$$
We consider the fibre product
$${\rm Spec}\big(A[1/t_1\ldots t_\ell]\big)\times_{\cMg} \cGMg .$$
This is finite \'etale over $A[1/t_1\ldots t_\ell]$ hence affine. The
normalization of this over ${\rm Spec}(A)$ is ${\rm Spec}(B)$; here $B$ is
a product of complete local rings finite over $A$, ramified only over
$t_1\ldots t_\ell=0$. Thus we have
$${\rm Spec}\big(B[1/t_1\ldots t_\ell]\big)=
{\rm Spec}\big(A[1/t_1\ldots t_\ell]\big)\times_{\cMg} \cGMg .$$

We claim that the morphism ${\rm Spec}(B)\to \bGMg$ identifies complete
local rings at the points of $\bGMg$ lying over $[C]\in \bMg(k)$. By
general theory we know that the completion of $\bMg$ at $[C]$ is
${\rm Spf}(A)/{\rm Aut}(C)$. A formal argument gives that the completion
of $\bGMg$ along $\varphi^{-1}([C])$ is isomorphic to
${\rm Spf}(B)/{\rm Aut}(C)$. Thus it suffices to show that
the action of ${\rm Aut}(C)$ on $\varphi^{-1}([C])$ is free. By comparing
levels, it suffices to show this for $G=(\msy Z/m\msy Z)^{2g}$; this is
the content of \cite[Proposition 3.5]{De}. See also references in remark
below.
%The result follows from the fact
%that an automorphism of $C$ is the identity if it acts trivially on
%$H^1(C, \msy Z/m\msy Z)$. See also \cite{Po,Mo,GO} and see remark below.

We use Abhyankar's lemma which asserts
that any normal local domain $A'$, finite generically
\'etale over $A$ and ramified only along $t_1\ldots t_\ell=0$ is contained
in $A[t_i^{1/n_i}]$ for some $n_1,\ldots,n_\ell$ relatively prime to $p$.
In addition, it is easily seen that $A'$ is formally smooth over $W(k)$
if and only if $A'$ is actually equal to $A[t_i^{1/n_i}]$
for some $n_1,\ldots,n_\ell$.

Let $C_K$ denote the lift of $C$ to $K=Q\big(W(k)\big)$ given by setting
$t_1=\ldots=t_{3g-3}=0$. The homomorphism
$A\to K[[t_1,\ldots,t_{3g-3}]]$ is such that
$$B\otimes_A K[[t_1,\ldots,t_{3g-3}]]$$
describes the complete local rings of $\bGMg$ along $\varphi^{-1}([C_K])$.
The result follows by comparing $B$ to this ring.\end{proof}

\begin{remark} {\rm The arguments above actually show that in the
situation of the proposition the stacks $\bcGMg$ are schemes, i.e.,
$\bcGMg\cong \bGMg$. Thus we get a stable curve over $\bGMg$ from
the morphism $\bGMg\cong \bcGMg\to \bcMg$; this also follows
from the existence of such a stable curve in the case
of abelian level structures, see \cite[Thm. 10.9]{Po},
 \cite[page 12]{GO} and \cite[Bemerkung 1]{Mo}.}\end{remark}

\section{Monodromy along the boundary}

In this section we study the moduli spaces $\bGMg$ along the boundary.
To do this it suffices to understand the monodromy along the boundary
on the relative fundamental group of the universal curve over $\cMgan$.

\subsection{The results}\label{results}

Let us formulate the main result. Let $\Pi$ denote the standard fundamental
group of a compact Riemann surface of genus $g$. We fix natural numbers
$k,n$ with $k\geq 1$ and $n\geq 3$. We will consider moduli of curves
of genus $g$ with a Teichm\"uller structure of level $G$, where
$$G=\Pi\big/ \Pi^{(k+1),n}.$$
(For notations see \ref{notation}.) By a result of Labute \cite{Labute},
the quotients $\Pi^{(k)}/\Pi^{(k+1)}$ are finitely generated free
abelian groups. Thus $G$ has a filtration whose successive quotients
are finite abelian groups of exponent $n$. Any prime dividing $\# G$
also divides $n$. Further, there is a surjection $G\to (\msy Z/n\msy Z)^{2g}$.
By Section 2 we get a moduli scheme
$$\bGMg\longrightarrow {\rm Spec}(\msy Z[1/n])$$
whose interior $\GMg$ classifies smooth genus $g$ curves $C$ with a surjection
$\pi_1(C)\to G$ given up to inner automorphisms.

\begin{theorem}\label{glad}
Suppose $k\in \{1,2,3\}$ and $n\geq 3$. The structural morphism
$\bGMg\to {\rm Spec}(\msy Z[1/n])$ is smooth if and only if
\begin{itemize}
\item{$k=1$} and $g=2$,
\item{$k=2$} and $n$ is odd,
\item{$k=3$} and $n$ is odd or $n$ is divisible by $4$.
\end{itemize}
Furthermore, if $k \geq 4$, $n \geq 3$ and $n$ relatively prime to $6$,
then $\bGMg\to {\rm Spec}(\msy Z[1/n])$ is smooth.
\end{theorem}

\begin{remark}{\rm The abelian case of this theorem, i.e., the case $k=1$,
has been proven by Mostafa \cite{Mo} and Van Geemen-Oort \cite{GO}.

A smooth and compact cover $\overline{M_g[_2^n]}^{an} \rightarrow \bMgan$
with $n$ even and $n\geq 6$
 has been constructed by Looijenga using Prym level structures,
see \cite{Looijenga}. This has been extended to higher level structures
$\overline{ M_g[_q^n]}^{an}$ by one of us (Pikaart). Using an
analog of Proposition \ref{criterium} this result may be extended
to an open set of ${\rm Spec}(\msy Z)$, at least if $q\geq 3$.

We note that $\bGMg\to\bMg[1/n]$ is a Galois cover\label{galois}
with group ${\rm Aut}(G)/{\rm Inn}(G)$ if $g\geq 3$ (if $g=2$ one has to divide
out an additional $\{\pm 1\}$). This follows since for
any surjection $\Pi\to G$ the kernel is equal to $\Pi^{(k),n}$. However,
$\bGMg\otimes\msy C$ need not be connected, see \cite[5.13]{DM} for
a description of the set of connected components.}\end{remark}

To prove the theorem, it suffices to consider $\bGMgan$, see Proposition
\ref{criterium}. We already know that $\GMgan$ is smooth. As in the
proof of Proposition \ref{criterium} we describe analytic neighbourhoods of
points in the boundary.

Let $C$ be a complex stable curve of genus $g$ with singular points
$\{ P_1,\ldots,P_\ell\}$. Let $\Gamma=\Gamma(C)$ be its dual graph;
an edge for each point $P_j$, a vertex for an irreducible component of $C$.
Let $ \pi: (\cC,C) \rightarrow (B,0)$ be a local
universal deformation of $C$,
where $B\subset \msy C^{3g-s}$ is a polydisc neighbourhood of $0$.
The coordinates $z_i$ are chosen such that $z_j=0$, $1\leq j\leq \ell$
parametrizes curves where the singular point $P_j$ subsists.
The discriminant locus $\Delta\subset B$ of $\pi$ is thus given
by $z_1\ldots z_\ell=0$. Put $U=B\setminus \Delta$, let $x\in U$
and choose $y\in \cC_x=\pi^{-1}(x)$.
The fundamental group of $U$ is an abelian group, freely generated by
simple loops
around the divisors $z_j=0$, thus naturally isomorphic to the free
abelian group on the edges of $\Gamma$, i.e., $\pi_1(U,x)
\cong \bigoplus_{e \in {\rm Edges}( \Gamma )} {\msy Z} e$.
The map $\cC|_U\to U$ is a locally trivial fibration, hence we have
the exact sequence
$$1 \longrightarrow \pi_1(\cC_x,y) \longrightarrow \pi_1(\cC|_U,y)
\longrightarrow \pi_1(U,x) \longrightarrow 1.$$
(Use that $\pi_2(U)=(0)$.) This short exact sequence provides us
with the monodromy representation
$$\rho: \pi_1(U,x) \longrightarrow {\rm Out}\big(\pi_1(\cC_x,y)\big).$$

The points $P_j$ determine non-trivial distinct isotopy classes
of circles on $\cC_x$, which have pairwise disjoint
representatives $c_j$. The fundamental group of $U$ is also naturally
isomorphic to the free abelian group on these circles, $\pi_1(U,x)
\cong \bigoplus_{i=1}^l {\msy Z} c_i$. Under this identification
we have that
$$\rho(c_i)=D_{c_i},$$
where $D_{c_i}$ is the exterior automorphism of $\pi_1(C_x)$ given by a Dehn
twist (also written $D_{c_j}$) around the circle $c_i$ (see
\cite{Dehn},\cite{Lamotke}).

We will describe of a neighbourhood of a point in $\bGMgan$
lying above $[C]$. Let $Z$ be the fibre product
$$Z=U\times_{\Mgan}\GMgan.$$
The normalization of $B$ in the function field of $Z$ is denoted $\bar Z$.
Note that $Z\to U$ is a finite topological covering space given by the set
$$S=\hbox{Hom-surj}\big(\pi_1(\cC_x,y),G\big)\big/\pi_1(\cC_x,y)$$
with $\pi_1(U,x)$-action defined via $\rho$. As in the proof of Proposition
\ref{criterium} there is an action of ${\rm Aut}(C)$ on $\bar Z$ and
$\bar Z/{\rm Aut}(C)$ defines a neighbourhood of
$\varphi^{-1}([C])\subset \bGMgan$. As in that proof we get that
$\bGMgan$ is smooth along $\varphi^{-1}([C])$ if and only if $\bar Z$
is smooth. (Here we use again that ${\rm Aut}(C)$ acts freely
on $\varphi^{-1}([C])$.) Finally there is the following
criterion: $\bar Z$ is smooth
if and only if for all $s\in S$ we have
$${\rm Stab}(s)=\bigoplus_{e \in {\rm Edges}(\Gamma)}n_e {\msy Z} e $$
for certain $n_e\in \msy Z$. Notice that if $m \in \msy Z $, then $mF^0 +m F^1$
($F^i$ as below) is of this form, but $2mF^0 +m F^1$ is not.

We remark that the arguments above go through for arbitrary finite
groups $G$, with a surjection onto $(\msy Z/n\msy Z)^{2g}$. In order
to describe the stabilizers in our more special situation we introduce the
following notation:
\begin{eqnarray*}
\Akn &= &{\rm Ker}\left( {\rm Aut}(\Pi) \to {\rm
Aut}(\Pi/\Pi^{(k+1),n})\right),\\
\Ikn &= &{\rm Ker}\left( {\rm Inn}(\Pi) \to {\rm
Inn}(\Pi/\Pi^{(k+1),n})\right),\\
\Okn &= &{\rm Ker}\left( {\rm Out}(\Pi) \to {\rm
Out}(\Pi/\Pi^{(k+1),n})\right).
\end{eqnarray*}
Oda et al.\ consider a variant with $n=0$.
By choosing an isomorphism $\pi_1(\cC_x,y)\cong\Pi$ we may view $\rho$
as a map into ${\rm Out}(\Pi)$. It is clear that
${\rm Stab}(s)=\rho^{-1}(\Okn)$ for any $s\in S$ (use Remark \ref{galois}).
Therefore, Theorem \ref{glad} follows from Theorem \ref{mono} below.

We will describe a decreasing
filtration $F^i$ on $\bigoplus_{e \in {\rm Edges}( \Gamma)}
{\msy Z} e$.
An edge $e$ such that $\Gamma \backslash e$ is disconnected is called a
{\em bridge}. A bridge $b$ is said to {\em bound a genus one curve}
 if one of the two
components of $C_x \backslash \{ $
the circle corresponding to $b \}$ has genus one.
A pair of distinct edges $\{e,f\}$ is called a {\em cut pair} if neither $e$
nor
$f$
is a bridge and $\Gamma \backslash \{e,f \}$ is disconnected.
A subset $E$ of the edges of $\Gamma$ is called a {\em maximal cut system}
if $E$ contains at least one cut pair, any two elements of $E$ form a cut pair
and no element of $E$ forms a cut pair with an element outside $E$.
Let $B$ be the set of bridges of $\Gamma$ and let $B_1$ be the subset of $B$
consisting of bridges which bound genus one curves. Let $\{E_i\}_{i \in I}$
denote the maximal cut systems.Set
$$D_i:=Ker(\bigoplus_{e \in E_i}
\msy Z e \stackrel{deg}{\rightarrow } \msy Z).$$
We define a decreasing filtration $F^i$ on $\bigoplus_{e \in {\rm
Edges}(\Gamma)}\msy Z  e$, as follows:
\begin{eqnarray*}
F^0&= & \bigoplus\nolimits_{e \in {\rm Edges}(\Gamma)}{\msy Z}  e,\\
F^1&= & \bigoplus\nolimits_{i\in I} D_i \bigoplus
\left(\bigoplus\nolimits_{b \in B} {\msy Z}b\right), \\
F^2&= & \bigoplus\nolimits_{b \in B} {\msy Z}b ,\\
F^3&= & (0).\end{eqnarray*}
Furthermore, set $F^2_1:=\bigoplus_{b \in B_1} {\msy Z}b$. This refines the
filtration into $F^0\supset F^1\supset F^2\supset F^2_1\supset F^3=(0)$.
Here is the main result of this article.

\begin{theorem}\label{mono} Notations as above. For $n,l \in \msy Z$, define
$n_l:=n/ {\rm gcd}(l,n)$.
\[ \begin{array}{cll}
1. & \mbox{If $k=1$, then} & \rho^{-1}(\Okn ) =nF^0 +F^1 \\
2. & \mbox{If $k=2$, then} &  \rho^{-1}(\Okn )=nF^0 +n_2F^1+F^2 \\
3. & \mbox{If $k=3$ and $2 ||n$, then}
&  \rho^{-1}(\Okn )=nF^0 + \frac{1}{2}nF^1+n_2F^2 +n_6F^2_1 \\
 & \mbox{If $k=3$ and $n$ is odd or $4 | n$, then}
& \rho^{-1}(\Okn ) =nF^0 + nF^1+n_2F^2 +n_6F^2_1\\
4. & \mbox{If $k \geq 4$ and $(n,6)=1$, then} & \rho^{-1}(\Okn ) =nF^0
\end{array} \] \end{theorem}

\begin{remark} {\rm The case $k=1$ has been proven by Brylinski
(\cite{Brylinski}). The case $k \geq 4$ follows from 3, the easy inclusions
of Section \ref{zeer easy} and the inclusions $\Okn  \subset O^{(l),n}$ if
$k \geq l$.

If we take $n=0$ then Theorem \ref{mono}
reduces to the non-pointed case of \cite[Main Theorem]{Oda}.} \end{remark}

\section{Description of Dehn twists and easy inclusions} \label{zeer easy}

In this section we will prove the inclusions ``$\supset$'' from Theorem
\ref{mono}.

Let $\Gamma$ be the graph of a stable curve and $(S,\{ c_i \})$ the smooth
model
with a set of circles as described in Subsection \ref{results}. We will
describe the Dehn twist associated to a bridge, cut pair or circle
separatedly.

\subsection{Bridges}\label{bridges}

Let $b$ be a bridge on $S$, let $g$ be the genus of $S$. Modulo a
homeomorphism the situation looks as follows:
\hfill \vspace{4 cm}
$$\hbox{\sl Fig.~1}$$
Cutting $S$ along $b$ yields the decomposition $S= S_1 \cup S_2$.
Let $g_i$ be the genus of $S_i$.
Choose a base point $p$ in $S_1$ and standard generators
$\alpha_{\pm i},$ $1 \leq i \leq g$, for $\pi :=\pi_1(S,p)$ such that
$\alpha_{\pm i}$ are in $S_1$ if $i \leq g_1$ and $\alpha_{\pm i}$
for $i\geq g_1+1$ hits $b$ exactly twice. We set
$v=[\alpha_1,\alpha_{-1}] \cdots [\alpha_{g_1},
\alpha_{-g_1}]$, it is freely homotopic to $b$ for a suitable orientation
of $b$.
We list the effect of the Dehn twist $D_b$ on the standard generators:
%(we set $v=[\alpha_1,\alpha_{-1}] \cdots [\alpha_{g_1},\alpha_{-(g_1)}]$)

\[ \begin{array}{ll}
i \in \{\pm 1,\ldots, \pm  g_1\} & {D}_{b}: \alpha_i \mapsto \alpha_i ,\\
i \in \{ \pm ( g_1 +1),\ldots,\pm g \} & {D}_{b}: \alpha_i \mapsto v^{-1}
\alpha_i v
\end{array} \]
This gives for the mth-powers:

\[ \begin{array}{ll}
i \in \{\pm 1,\ldots, \pm  g_1\} & {D}^m_{b}: \alpha_i \mapsto \alpha_i ,\\
i \in \{ \pm ( g_1 +1),\ldots,\pm g \} & {D}^m_{b}: \alpha_i \mapsto v^{-m}
\alpha_i v^m
\end{array} \]

We see that $D_b(\alpha_{\pm i})\alpha_{\pm i}^{-1} \in \pi^{(3)}$
for all $i$. This proves that $\rho(b) \in O^{(2)}$ and thus
by linearity $\rho(F^2) \subset O^{(2)}$.

Suppose $b$ bounds a genus one curve, say $S_2$ has genus one.
In this case $v = [\alpha_{-g},\alpha_g]$ and we only have to consider
$$D_b^m(\alpha_{\pm g})\alpha_{\pm g}^{-1}=
[[\alpha_g,\alpha_{-g}]^m,\alpha_{\pm g}]\equiv
[[\alpha_g,\alpha_{-g}],\alpha_{\pm g}]^m \bmod \pi^{(4)}.$$
To prove this element lies in $\pi^{(4),n}$ we define $f:F(x,y,z) \ra \pi$ by
$x \mapsto \alpha_g$, $y \mapsto \alpha_{-g}$ and $z \mapsto 1$.
Then $f([[x,y],x])=[[\alpha_g,\alpha_{-g}],\alpha_{ g}]$
and $f([[x,y],y])=[[\alpha_g,\alpha_{-g}],\alpha_{ -g}]$.
{}From Lemma \ref{berekening} we see that $n_6|m$ implies
that the element above lies in $\pi^{(4),n}$.
Thus $\rho(n_6F^2_1) \subset O^{(3),n}$.

Suppose $b$ does not bound a genus one curve. In that case we have
\begin{eqnarray*}
D_b^m(\alpha_{\pm i})\alpha_{\pm i}^{-1}&=& [v^m,\alpha_{\pm i}^{-1}]\ \equiv\
[v,\alpha_{\pm i}]^{-m}\\
&\equiv & [[\alpha_{1},\alpha_{-1}],\alpha_{\pm i}]^{-m}\cdots
[[\alpha_{g_1},\alpha_{-g_1}],\alpha_{\pm i}]^{-m}
\end{eqnarray*}
for $i>g_1$ in $\pi/\pi^{(4)}$. Define $g_j:F(x,y,z) \ra \pi$ by $x \mapsto
\alpha_j,~y \mapsto \alpha_{-j}$
 and $z \mapsto \alpha_{\pm i}$. Then $g_j([[x,y],z]^m)=
[[\alpha_j,\alpha_{-j}],\alpha_{\pm i}]^m$.
{}From Lemma \ref{berekening} we see that
$n_2 |m$ implies $\rho(mb) \in O^{(3),n}$, and
thus $\rho(n_2F^2) \subset O^{(3),n}$.

\subsection{Edges which are not bridges}\label{not bridges}

Let $c$ be a circle on $S$ which is not a bridge. Modulo a homeomorphism the
situation looks as follows:
\hfill \vspace{4 cm}
$$\hbox{\sl Fig.~2}$$
We choose a point $p$ in $S$ and standard generators $\alpha_{\pm i}$ for
$\pi =\pi_1(S,p)$ such that $\alpha_{-g}$ is the only one intersecting $c$
and $\alpha_g$ is freely homotopic to $c$.
The action of the Dehn twist becomes:
$D_c(\alpha_i)=\alpha_i$ if $i \neq -g$ and $D_c(\alpha_{-g}) =\alpha_g
\alpha_{-g}$.
Thus $D_c^m(\alpha_i)\alpha_i^{-1}=1$ or $\alpha_g^m$.
Evidently, if $n|m$ then $\rho(mc) \in O^{(k),n}$ for all $k$.
Together with the results of \ref{bridges} this gives
$\rho(nF^0)\subset O^{(k),n}$ for all $k$.

\subsection{Cut systems}\label{cut systems}

Let $e_1,e_2$ be a cut pair on $S$. Modulo a homeomorphism the situation is
as follows:
%ruimte voor plaatje
\hfill \vspace{4 cm}
$$\hbox{\sl Fig.~3}$$
Cutting $S$ along $e_1$ and $e_2$ yields the decomposition $S=S_1 \cup
S_2$. Let $g_1$ be the genus of $S_1$. Choose a base point $p$ in $S_1$.
We take standard generators $\alpha_{\pm i}$ for $\pi_1(S,p)$, such that
 for $i \in \{ \pm 1,\ldots, \pm g_1, g_1+1\} $, $ \alpha_i $ is in $S_1$,
for $i \in \{ \pm (g_1+2),\ldots, \pm g \} $, $ \alpha_i $ enters $S_2$ via
$e_1$
and leaves $S_2$ also via $e_1$  and
$\alpha_{-(g_1+1)}$ enters $S_2$ via $e_1$ and leaves $S_2$ via $e_2$.
Furthermore, we want $\alpha_{g_1+1}$ to be freely homotopic to $e_2$, for
 a suitable orientation of $e_2$. It follows that after suitable orientation
 of $e_1$ the loop $[\alpha_1,\alpha_{-1}] \cdots
[\alpha_{g_1},\alpha_{-g_1}]\alpha_{g_1+1}^{-1}$ is freely homotopic to
${e_1}$. Let us write $v=[\alpha_1,\alpha_{-1}] \cdots
[\alpha_{g_1},\alpha_{-g_1}]$.
We list the effect of the Dehn twists ${D}_{e_i}$ on these generators:
\[ \begin{array}{ll}
\mbox{$i \in \{ \pm 1,\ldots, \pm g_1, g_1+1\} $}&
	  \left\{ \begin{array}{l} {D}_{e_1}: \alpha_i \mapsto \alpha_i \\
		{D}_{e_2}: \alpha_i \mapsto \alpha_i
		    \end{array} \right. \\[12pt]
\mbox{$i = -(g_1+1)$} &
	\left\{ \begin{array}{l} {D}_{e_1}: \alpha_i \mapsto
\alpha_i v \alpha_{g_1+1}^{-1}\\
	   {D}_{e_2}: \alpha_i \mapsto
 \alpha_{g_1+1}^{-1}\alpha_i
		\end{array} \right. \\[12pt]
\mbox{$i \in \{ \pm (g_1+2),\ldots, \pm g \}$} &
\left\{ \begin{array}{l} {D}_{e_1}: \alpha_i \mapsto \alpha_{g_1+1} v^{-1}
\alpha_i v \alpha_{g_1+1}^{-1}\\
	       {D}_{e_2}: \alpha_i \mapsto \alpha_i
	       \end{array} \right.
\end{array} \]
Thus we get the following formulae for $D_{e_2}D_{e_1}^{-1}$:
\[ \begin{array}{ll}
|i| \leq g_1  & D_{e_2}D_{e_1}^{-1}(\alpha_i)=\alpha_i, \\[2pt]

i=g_1+1 & D_{e_2}D_{e_1}^{-1}(\alpha_i)=\alpha_i, \\[2pt]

i = -(g_1+1)& D_{e_2}D_{e_1}^{-1}(\alpha_i)
 =\alpha_{g_1+1}^{-1}\alpha_i\alpha_{g_1+1}v^{-1}\alpha_i^{-1}\alpha_i
=[\alpha_{g_1+1},\alpha_i^{-1}]v^{-1}[v^{-1},\alpha_i^{-1}]\alpha_i,\\[2pt]

|i| \geq g_1+2 & D_{e_2}D_{e_1}^{-1}(\alpha_i)
=v\alpha_{g_1+1}^{-1}\alpha_i\alpha_{g_1+1}v^{-1}\alpha_i^{-1}\alpha_i
=v[\alpha_{g_1+1},\alpha_i^{-1}]v^{-1}[v^{-1},\alpha^{-1}]\alpha_i
\end{array} \]
and for the mth powers:
\[ \begin{array}{ll}
i \in \{ \pm 1,\ldots, \pm g_1, g_1+1\}&
(D_{e_2}D_{e_1}^{-1})^m(\alpha_i)=\alpha_i,\\[2pt]

i = -(g_1+1)& (D_{e_2}D_{e_1}^{-1})^m(\alpha_i)
 =\alpha_{g_1+1}^{-m}\alpha_i (\alpha_{g_1+1}v^{-1})^m,\\[2pt]

i \in \{ \pm (g_1+2),\ldots, \pm g \}&
(D_{e_2}D_{e_1}^{-1})^m(\alpha_i)
=(v\alpha_{g_1+1}^{-1})^m \alpha_i (\alpha_{g_1+1}v^{-1})^m.
\end{array} \]

This proves that for a cut pair  $\{e_1,e_2\}$ we have
$\rho(e_2-e_1) \in O^{(1)}$. It follows from this and \ref{bridges} that
$\rho(F^1) \subset O^{(1)}$, because
$F^1$ is generated by elements of the form $e_2-e_1$ for a cut pair
$\{e_1,e_2\}$ and the elements $b,~b \in B$.

To finish the proof of the inclusions ``$\supset$''
we have to show that $\rho(me_2-me_1)\in O^{(2),n}$ if
$n_2|m$ and $\rho(me_2-me_1)\in O^{(3),n}$ if $n|m$, or
$\frac{1}{2}n|m$ in case $2||m$. We have to show that the
divisibility conditions imply
$(D_{e_2}D_{e_1}^{-1})^m(\alpha_i)\alpha_i^{-1}\in \pi^{(3),n}$,
respectively $\pi^{(4),n}$. We will make computations in $\pi$
modulo $\pi^{(4)}$. In the case that $i\in\{\pm1,\ldots,\pm g_1, g_1+1\}$
there is nothing to prove. We have
$(\alpha_{g_1+1}v^{-1})^m\equiv
\alpha_{g_1+1}^mv^{-m}[v^{-1},\alpha_{g_1+1}]^{\frac{1}{2}m(m-1)}
$ mod $ \pi^{(4)}$, as one can prove by induction.
{}From equality 5 of Subsection \ref{group structure} we have
$$\alpha_{g_1+1}^{-m}\alpha_i\alpha_{g_1+1}^m \alpha_i^{-1}\equiv
 [\alpha_{g_1+1},\alpha_i^{-1}]^m
 [[\alpha_{g_1+1},\alpha_i^{-1}],\alpha_{g_1+1}]^{\frac{1}{2}m(m-1)}.$$
Thus for $i=-(g_1+1)$ we get
\begin{eqnarray*}
&&(D_{e_2}D_{e_1}^{-1})^m(\alpha_i)\alpha_i^{-1}\\
&\equiv&\alpha_{g_1+1}^{-m} \alpha_i (\alpha_{g_1+1}v^{-1})^m \alpha_i^{-1}\\
&\equiv&\alpha_{g_1+1}^{-m} \alpha_i
\alpha_{g_1+1}^mv^{-m} [v^{-1},\alpha_{g_1+1}]^{\frac{1}{2}m(m-1)}
 \alpha_i^{-1}\\
&\equiv&{[\alpha_{g_1+1},\alpha_i^{-1}]}^{m}
[[\alpha_{g_1+1},\alpha_i^{-1}],\alpha_{g_1+1}]^{\frac{1}{2}m(m-1)}
v^{-m} [v, \alpha_i]^m [v,\alpha_{g_1+1}]^{-\frac{1}{2}m(m-1)}\\
&\equiv&
{[\alpha_{g_1+1},\alpha_i^{-1}]}^{m}
[[\alpha_{g_1+1},\alpha_i^{-1}],\alpha_{g_1+1}]^{\frac{1}{2}m(m-1)}
[\alpha_1,\alpha_{-1}]^m \cdots [\alpha_{g_1},\alpha_{-g_1}]^m
 [[\alpha_1,\alpha_{-1}],\alpha_i]^m \cdots \\
&&
{[[\alpha_{g_1},\alpha_{-g_1}],\alpha_i]}^m
{[[\alpha_1,\alpha_{-1}],\alpha_{g_1+1}]}^{-\frac{1}{2}m(m-1)} \cdots
[[\alpha_{g_1},\alpha_{-g_1}],\alpha_{g_1+1}]^{-\frac{1}{2}m(m-1)}.
\end{eqnarray*}
In this product there are, apart from factors in $\pi^{(3)}$,
a number of terms of the form
$[x,y]^m$, $x,y\in \pi$. By the divisibility conditions
all these are in $\pi^{(4),n}$, since this is true in the
case of a free group $\langle x, y, z\rangle$ by Lemma \ref{berekening}.
Thus it is clear that the whole product lies in $\pi^{(3),n}$ for
$n_2|m$.

The other terms in the product are of the form
$[[x,y],z]^{\frac{1}{2}m(m-1)}$, where $x,y,z\in \pi$.
It follows from Lemma \ref{berekening}, by mapping $G=\langle x, y, z\rangle$
into $\pi$ that these terms are in $\pi^{(4),n}$ as soon as
$n_2| \frac{1}{2}m(m-1)$. This is equivalent to the condition
$n|m$ or $\frac{1}{2}n|m$ in case $2||n$.

For $i \geq g_1+2$ we have:
\begin{eqnarray*}
(D_{e_2}D_{e_1}^{-1})^m(\alpha_i)\alpha_i^{-1}&=&
(v\alpha_{g_1+1}^{-1})^m \alpha_i (\alpha_{g_1+1}v^{-1})^m
\alpha_i^{-1} \\
&\equiv&
v^m \alpha_{g_1+1}^{-m} [\alpha_{g_1+1}^{-1},v]^{\frac{1}{2}m(m-1)}\alpha_i
\alpha_{g_1+1}^mv^{-m} [v^{-1},\alpha_{g_1+1}]^{\frac{1}{2}m(m-1)}
 \alpha_i^{-1}\\
&\equiv &
{[\alpha_{g_1+1},\alpha_i^{-1}]}^m
[[\alpha_{g_1+1},\alpha_i^{-1}],\alpha_{g_1+1}]^{\frac{1}{2}m(m-1)}
\end{eqnarray*}
The same argument works to show that our divisibility conditions
imply this product lies in either $\pi^{(3),n}$ or $\pi^{(4),n}$.

\section{Completion of proof}

In this section we will use the following argument.
Suppose we are given a morphism $f:S \ra S'$ of compact connected oriented
surfaces and a circle $c$ on $S$. Assume that $f$ maps a tubular neighbourhood
of $c$ isomorphically into $S'$, denote by $c'$ the image of $c$.
Then we have: $f_\ast \circ D_c =D_{c'} \circ f_\ast$.

Let $\Gamma$ be the graph of a stable curve and $(S,\{ c_i \})$
the smooth model with a set of circles as described in Subsection
\ref{results}.
Let $H$ be the set of circles which are not bridges and are not
involved in any cut system.
Let $B$ be the set of bridges of $\Gamma$ and let $\{ E_i|i=1,\ldots,s \}$
be the set of maximal cut systems.
Choose a numbering $E_i= \{ e_{i,0},\ldots,e_{i,f_i} \}$ for each $i$.
Remark that ${\rm Edges}(\Gamma)=H\cup \bigcup E_i\cup B$.

\subsection{The case $k=1$}

We want to prove the inclusion ``$\subset$'' for $k=1$.
Suppose $\sigma \in \oplus {\msy Z}e$ is a counterexample which involves a
minimal number of edges. This means that $\rho(\sigma)$ lies in $O^{(1),n}$
but not all coefficients of $\sigma$ are divisible by $n$; the number
of nonzero coefficients is minimal. We write $\sigma =\sum n_c c$.
We may subtract elements of $nF_0+F_1$ from $\sigma$:
we already know that $nF_0+F_1$ maps into $O^{(1),n}$.
Thus, by minimality, we may suppose that:
(\romannumeral1) none of the nonzero coefficients $n_c$ is divisible by $n$,
(\romannumeral2) none of the circles $c$ is a bridge, and
(\romannumeral3) of each cutsystem $E_i$ at most one element occurs in
$\sigma$.  Modulo a homeomorphism the situation looks as follows:
\hfill \vspace{4 cm}
$$\hbox{\sl Fig.~4}$$
Take a point $p$ near $c$. It is clear that we can find a loop $\alpha$
which intersects $c$ exactly once and none of the other edges involved in
$\sigma$. But now we have a contradiction, because by Section \ref{not bridges}
we know: $\rho(\sigma)(\alpha)\alpha^{-1}= c^{n_c} \notin \pi^n$
if $n$ does not divide $n_c$.
Thus we have proven $\rho^{-1}(O^{(1),n})\subset nF^0+F^1$, equality follows
from sections \ref{not bridges} and \ref{cut systems}.

\subsection{The case $k=2$}\label{k is 2}

We want to prove the inclusion ``$\subset$'' for $k=2$.
Let $\sigma=\sum_{b \in B}m_b b
+\sum_{i,j}m_{i,j}(e_{i,j}-e_{i,0})
+\sum_i m_i e_{i,0} + \sum_{c \in H} n_c c$ be such that $\rho(\sigma ) \in
O^{(2),n}$. By the above we know that $n|m_i$ and $n|n_c, c\in H$. We have
to show that $n_2|m_{i,j}$ for all possible $i,j$.
Suppose this does not hold and suppose furthermore that $\sigma$ is
a minimal counterexample with respect to the number of edges involved.
Arguing as above we may assume that
$\sigma =\sum_{i,j}m_{i,j}(e_{i,j}-e_{i,0})$, where no (nonzero)
$m_{i,j}$ is divisible by $n_2$.

\begin{proposition} \label{nice cut system}
There is a maximal cut system which is such that if we cut $S$ along
every element of the cut system, one connected component contains all other
cut systems. \end{proposition}
\begin{proof} (Cf.\ \cite{Oda}, Lemma 5.3.)
Take any cut system $E$. If it does not have the required
property, that means that there is more than one  connected component of $S
\setminus E$, say $S_{E,1}, \dots , S_{E,l}$, which contain maximal
 cut systems. Suppose $S_{E,j}$ contains $a_j$ maximal cut systems and
suppose $a_1= \mbox{ min}_j \{ a_j \}$. Take a maximal cut system $F$ in
$S_{E,1}$. Clearly there is a connected component of $S\setminus F$ which
contains at least $\sum_{j=2}^la_j+1$ components, namely the one containing
$E$. This component contains more maximal cut systems then the one we started
with.
Continuing in this way, we arrive at our result.
\end{proof}

Let $E_1=\{e_{1,0}, \dots, e_{1,f_1} \}$ be a
maximal cut system such that one of the two components bounded by $e_{1,0}$ and
$e_{1,f_1}$, contains all other cut systems involved in $\sigma$.
We suppose these edges are numbered cyclically: one of the connected
components of $S\setminus \{e_{1,i},e_{1,i+1}\}$ contains all other
$e_{1,j}$ (see figure 6). If $E_1$ contains more than two edges,
we proceed as follows. Replace
the component of $S \setminus \{ e_{1,1},~e_{1,f_1} \}$
which does not contain $e_{1,0}$
by a cylinder to get an oriented surface $S'$.
There is a continuous map $f: S\to S'$ such that for all elements
$e_{1,j}$, except $e_{1,0}$, we have that $f(e_{1,j})$ is homotopic
to $f(e_{1,1})$. As explained above, we get that
$f_\ast \circ \rho(\sigma)=\rho(f_\ast\sigma)\circ f_\ast$.
Thus $\rho(\sigma)\in O^{(2),n}$
implies $\rho(f_\ast \sigma)\in O^{(2),n}$.
It is also
clear that $f_\ast(\sigma)= m_{1,1}(f(e_{1,1})-f(e_{1,0}))+
\sum_{i>1} m_{i,j}(f(e_{i,j})-f(e_{i,0}))$. We are reduced to the case
that $E_1=\{e_{1,0},e_{1,1}\}$ consists of two elements.

The situation now is as in \ref{cut systems} (figure 3) where
all the cut systems $E_i$, $i>1$, lie in the component $S_2$.
The generator $\alpha_{-(g_1+1)}$ may be chosen such that it does not
intersect the circles $e_{i,j}$, $i>1$. Thus we get from the formulae
of Section \ref{cut systems} that
\begin{eqnarray*}
\rho(\sigma)(\alpha_{-(g_1+1)})\cdot \alpha_{-(g_1+1)}^{-1}&
=&(D_{e_{1,1}}D_{e_{1,0}}^{-1})^{m_{1,1}}(\alpha_{-(g_1+1)})
\cdot \alpha_{-(g_1+1)}^{-1}\\
&=&\alpha_{g_1+1}^{-m_{1,1}}\alpha_{-(g_1+1)}(\alpha_{g_1+1}v^{-1})^{m_{1,1}}
\alpha_{-(g_1+1)}^{-1}.\end{eqnarray*}
We define $f:\pi\to G=\langle x,y,z\rangle$ by
$\alpha_1, \alpha_{-g} \mapsto x$
, $\alpha_{-1}, \alpha_g \mapsto y$, $\alpha_{g_1+1} \mapsto z$,
other generators $ \mapsto 1$. (The defining relation for
$\pi$ is indeed mapped to 1.) The expression above is mapped
to $z^{-m_{1,1}}\cdot(z\cdot [y,x])^{m_{1,1}}$.
Modulo $G^{(3)}$ this equals $[y,x]^{m_{1,1}}$, hence
$\rho(\sigma)\in O^{(3),n}$ implies $n_2|m_{1,1}$ (see Corollary
\ref{free comp}). In this way we prove the desired contradiction;
the inclusion ``$\subset$'' for $k=2$ follows.

\subsection{The case $k=3$}

We want to prove the inclusion ``$\subset$'' for $k=3$.
Let $\sigma=\sum_{b \in B}m_b b
+\sum_{i,j}m_{i,j}(e_{i,j}-e_{i,0})
+\sum_i m_i e_{i,0} + \sum_{c \in H} n_c c$ be such that
$\rho(\sigma) \in O^{(3),n}$. We have to show that
$n_2|m_{i,j}$ or $n|m_{i,j}$, for all $i,j$,
depending on $n$ being exactly divisible by $2$ or not;
and $n_6$ respectively $n_2$ divides $m_b$ depending on $b \in B_1$ or not.
Suppose this does not hold and suppose furthermore that $\sigma$ is
a minimal counterexample with respect to the number of edges involved.
Thus $\sigma =\sum_{b \in B}m_b b +\sum_{i,j}m_{i,j}(e_{i,j}-e_{i,0})$
and none of the nonzero $m_b$, $m_{i,j}$ satisfy the divisibility conditions.

Case 1. For some $b$ involved in $\sigma$ one of the connected
components, say $S'$, of $S\setminus \{b\}$ contains no other
edges involved in $\sigma$. The situation looks as follows.
\hfill \vspace{4.8 cm}
$$\hbox{\sl Fig.~5}$$
We take a basepoint $p$ in $S'$ and standard generators
$\alpha_{\pm i}$ which are loops in $S'$ for $i \leq g(S')$
and which are not loops in $S'$ for $i > g(S')$.
We may choose $\alpha_{g(S')+1}$ such that it does not intersect any
circles involved in $\sigma$ but $b$. To see this let $S''$ denote
the connected component of $S\setminus \{c\in \sigma, c\not = b\}$
containing $S'$. If $g(S'')>g(S')$, then there is a `hole' between
$b$ and the boundary of $S''$, which we can take to be hole number
$g(S')+1$ and take $\alpha_{g(S')+1}$ accordingly. This is the
case for example if $S''=S$ (i.e., if $\sigma = m_b b$)
or if $S''$ is bounded by another bridge $b'$. If $g(S'')=g(S')$,
then $S''$ is bounded by a cut pair $\{e_1, e_2\}$ (belonging to
some maximal cut system $E_i$), the part $S''\setminus S'$
looks like a pair of pants. In this case we can choose $\alpha_{g(S')+1}$
to go around a pants leg of $S''\setminus S'$, for example freely homotopic
to $e_2$ as in Section \ref{cut systems}.

It is now clear that $\rho(\sigma)(\alpha_{g(S')+1})=
D_b^{m_b}(\alpha_{g(S')+1})$ and by Section \ref{bridges} we get
$$ \xi:= \rho(\sigma)(\alpha_{g(S')+1})\cdot \alpha_{g(S')+1}^{-1}=
v^{-m_b}\alpha_{g(S')+1}v^{m_b}\alpha_{g(S')+1}^{-1}$$
with $v=[\alpha_1,\alpha_{-1}]\cdot\ldots\cdot
[\alpha_{g(S')},\alpha_{-g(S')}]$.
Thus $\rho(\sigma)\in O^{(3),n}$ implies $\xi\in \pi^{(4),n}$.

If $b$ bounds a curve of genus $1$, say $g(S')=1$ (the case $g-g(S')=1$
is similar), then $v=[\alpha_1,\alpha_{-1}]$. We define $f :\pi \to
G=\langle x,y,z\rangle$ by $\alpha_1,\alpha_{-2}\mapsto x$,
$\alpha_{-1},\alpha_{2}\mapsto y$ and other generators $\mapsto 1$.
We see that $f(\xi)= [x,y]^{-m_b}y[x,y]^{m_b}y^{-1}\equiv [[x,y],y]^{-m_b}$
modulo $G^{(4)}$. Thus $\xi\in \pi^{(4),n}$ implies $n_6 | m_b$ by
Lemma \ref{berekening}.

If $b$ does not bound a curve of genus 1 (i.e., $g(S')\geq 2$ and
$g-g(S')\geq 2$) we define $g:\pi \to G=\langle x,y,z\rangle$ by
$\alpha_1, \alpha_{-g} \mapsto x$, $\alpha_{-1}, \alpha_g \mapsto y$
, $\alpha_{g(S')+1} \mapsto z^{-1}$ and other generators $\mapsto 1$.
We get $g(\xi)=[x,y]^{-m_b}z^{-1}[x,y]^{m_b}z\equiv [[x,y],z]^{m_b}$.
Application of Lemma \ref{berekening} gives $n_2|m_b$ as desired.

Case 2: case 1 does not occur. This means that every
connected component of $S$ cut out by one bridge involved in $\sigma$ contains
a maximal cut system. Choose a maximal cut system as in Proposition \ref{nice
cut system}. Let $E_1=\{e_{1,0}, \dots, e_{1,f_1} \}$
be this cut system. Suppose the edges are numbered cyclically and
suppose the connected component cut out by $e_{1,0}$ and $e_{1,1}$
which does not contain any edge of $E_1$ contains all other
maximal cut systems. Call this component $S'$. By the choice of $E_1$ and
the assumption,  we get that $S \setminus S'$ does not contain any bridges
involved in $\sigma$. The situation thus looks as shown below:
\vspace{6.5 cm}
$$\hbox{\sl Fig.~6}$$
Arguing as in the proof for $k=2$,
we arrive at a curve with a cut pair $e_{1,0}$ and $e_{1,1}$ which is a maximal
cut system with the property that one of the two components cut out
by the cut pair contains all other edges involved in the counter example
$\sigma$. Considering $f : \pi \to G$ as in
Section \ref{k is 2}, we get, as in case 1, that the element
$z^{-m_{1,1}}\cdot (z\cdot [y,x])^{m_{1,1}}$ lies in $G^{(4),n}$.
Note that
$$z^{-m_{1,1}}\cdot (z\cdot [y,x])^{m_{1,1}}\equiv
[y,x]^{m_{1,1}}[[y,x],z]^{\frac{1}{2}m_{1,1}(m_{1,1}-1)}.$$
Thus we get our divisibility condition on $m_{1,1}$ by Lemma \ref{berekening}.

\section{Computations for a free group on three generators}

\label{group structure}
\label{berekeningen}

\begin{lemmatje} Let $G$ be a free group on three generators, $G=\langle
x,y,z\rangle$. Then
$$ G/G^{(4)} \cong \{(i_1,\ldots,i_{14})\in {\msy Z}^{14} \}$$
where multiplication on the right-hand side is given as follows:
$$(i_1,i_2,i_3,i_4,i_5,i_6,i_7,i_8,i_9,i_{10},i_{11},i_{12},i_{13},i_{14})
\cdot $$
$$(j_1,j_2,j_3,j_4,j_5,j_6,j_7,j_8,j_9,j_{10},j_{11},j_{12},j_{13},j_{14})=$$
\[ \left( \begin{array}{c}
i_1 +j_1 \\
 i_2 +j_2 \\
 i_3 +j_3 \\
i_4+j_4 -i_2 j_1 \\
 i_5+j_5-i_3 j_2 \\
i_6+j_6-i_3j_1  \\
i_7+j_7 +i_4j_1- \frac{1}{2}(j_1-1)j_1i_2 \\
i_8+j_8 +i_4j_2- \frac{1}{2}(i_2-1)i_2j_1 -i_2j_1j_2 \\
i_9+j_9 +i_4j_3+i_6j_2 -i_2i_3j_1-i_2j_1j_3-i_3j_1j_2 \\
i_{10}+j_{10}+i_5j_1+i_6j_2 -i_3j_1j_2 \\
i_{11}+j_{11}+i_5j_2 -\frac{1}{2}(j_2-1)j_2i_3 \\
i_{12}+j_{12}+i_5j_3-\frac{1}{2}(i_3-1)i_3j_2-i_3j_2j_3 \\
i_{13}+j_{13}+i_6j_1-\frac{1}{2}(j_1-1)j_1i_3 \\
i_{14}+j_{14}+i_6j_3-\frac{1}{2}(i_3-1)i_3j_1-i_3j_1j_3
\end{array} \right) \]
\end{lemmatje}

\begin{proof} We have the exact sequences
$$1 \rightarrow \Gk/\Gkeen \rightarrow G/ \Gkeen
\rightarrow G/ \Gk \rightarrow 1 ,$$
where $\Gk/\Gkeen$ is a finitely generated free abelian group
(\cite[Th.5.12]{Magnus}) of rank $N_i:=\frac{1}{i} \sum_{d|i} \mu(d) 3^{i/d}$
(\cite[Th.5.11]{Magnus}), where $\mu$ is the M\"obius-function.
Explicitly, we have, with $r:=[x,y]$, $s:=[y,z]$ and $t:=[x,z]$:
\[ \begin{array}{l}
G^{(1)}/G^{(2)} = \langle x,y,z\rangle_{ab} ~,\\
G^{(2)}/G^{(3)} = \langle r,s,t\rangle_{ab}~, \\
G^{(3)}/G^{(4)} =
\langle [r,x],[r,y],[r,z],[s,x],[s,y],[s,z],[t,x],[t,z]\rangle_{ab}~,
\end{array} \]
by the Jacobi-relation $[t,y] \equiv [r,z][s,x]$ mod $G^{(4)}$.
Now recall that if $a \in \Gk$ and $b \in G^{(l)}$ then $ab=ba[a,b]$
 with $[a,b] \in G^{(k+l)}$, so that
\begin{eqnarray} ab\equiv ba \mbox{  mod } G^{(k+l)}. \end{eqnarray}
Furthermore we have the following identities: (\cite[Th.5.1]{Magnus})
\begin{eqnarray}
\ [a,b]  & = & [b,a]^{-1},\\
\ [a,bc] & = & [a,c]~[a,b]~[[a,b],c], \\
\ [ab,c] & = & [a,c]~[[a,c],b]~[b,c].
\end{eqnarray}
It is clear now that any element $a$ of $G/G^{(4)}$ can be written in the form
$a=x^{i_{1}}\cdots [t,z]^{i_{14}}$ in a unique manner. Note
that the last eleven factors commute by (1).
We define a map $\phi:G/G^{(4)} \rightarrow \{(i_1,\ldots,i_{14}) \in {\msy
Z}^{14} \}$
by
$$x^{i_1}y^{i_2}z^{i_3}r^{i_4}s^{i_5}t^{i_6}[r,x]^{i_7}[r,y]^{i_8}[r,z]^{i_9}
[s,x]^{i_{10}}[s,y]^{i_{11}}[s,z]^{i_{12}}[t,x]^{i_{13}}[t,z]^{i_{14}}
\mapsto (i_1,\ldots,i_{14}) .$$
Before we start the computation, notice that modulo $G^{(4)}$ we have
\begin{eqnarray}
[a^i,b^j]\equiv
[a,b]^{ij}[[a,b],a]^{\frac{1}{2}ij(i-1)}[[a,b],b]^{\frac{1}{2}ij(j-1)},
\end{eqnarray}
as one proves easily by induction. We set $(n,m):=\frac{1}{2}(n-1)nm$.
The following identities hold in the group $G/G^{(4)}$:
\begin{eqnarray*}
&&x^{i_1}y^{i_2}z^{i_3}x^{j_1}y^{j_2}z^{j_3}\\
&=&
x^{i_1+j_1}x^{-j_1}y^{i_2}z^{i_3}x^{j_1}z^{-i_3}y^{-i_2}y^{i_2}z^{i_3}
y^{j_2}z^{j_3}\\
&=&
x^{i_1+j_1}[x^{j_1},z^{-i_3}y^{-i_2}]y^{i_2+j_2}y^{-j_2}z^{i_3}
y^{j_2}z^{-i_3}z^{i_3+j_3}\\
&=&
x^{i_1+j_1}[x^{j_1},y^{-i_2}][x^{j_1},z^{-i_3}][[x^{j_1},z^{-i_3}],y^{-i_2}]
y^{i_2+j_2}[y^{j_2},z^{-i_3}]z^{i_3+j_3}\\
&=&
x^{i_1+j_1}r^{-i_2j_1}[r,x]^{(j_1,-i_2)}[r,y]^{(-i_2,j_1)}t^{-i_3j_1}
[t,x]^{(j_1,-i_3)}[t,z]^{(-i_3,j_1)} \\
&&\ \ \ y^{i_2+j_2}s^{-i_3j_2}[s,y]^{(j_2,-i_3)}
[s,z]^{(-i_3,j_2)}z^{i_3+j_3}[t,y]^{i_2i_3j_1}\\
&=&
x^{i_1+j_1}y^{i_2+j_2}z^{i_3+j_3}r^{-i_2j_1}s^{-i_3j_2}t^{-i_3j_1}
[r,y]^{-i_2j_1(i_2+j_2)}[t,y]^{-i_3j_1(i_2+j_2)}[r,z]^{-i_2j_1(i_3+j_3)} \\
&&\ \ \ [t,z]^{-i_3j_1(i_3+j_3)}[s,z]^{-i_3j_2(i_3+j_3)}
[r,x]^{(j_1,-i_2)}[r,y]^{(-i_2,j_1)}[t,x]^{(j_1,-i_3)}[t,z]^{(-i_3,j_1)} \\
&&\ \ \
[s,y]^{(j_2,-i_3)}[s,z]^{(-i_3,j_2)}[r,z]^{j_1i_2i_3}[s,x]^{j_1i_2i_3}\\
&=&x^{i_1+j_1}y^{i_2+j_2}z^{i_3+j_3}r^{-i_2j_1}s^{-i_3j_2}t^{-i_3j_1}
[r,x]^{(j_1,-i_2)}[r,y]^{-(i_2,j_1)-i_2j_1j_2}[r,z]^{-i_2j_1j_3-i_3j_1j_2
-i_2i_3j_1} \\
&&\ \ \ [s,x]^{-i_3j_1j_2}[s,y]^{(j_2,-i_3)}[s,z]^{-(i_3,j_2)-i_3j_2j_3}
[t,x]^{(j_1,-i_3)}[t,z]^{-(i_3,j_1)-i_3j_1j_3}.
\end{eqnarray*}
Also
\begin{eqnarray*}
r^{i_4}s^{i_5}t^{i_6}x^{j_1}y^{j_2}z^{j_3}r^{j_4}s^{j_5}t^{j_6}&=&
x^{j_1}y^{j_2}z^{j_3}r^{i_4}s^{i_5}t^{i_6}[r^{i_4}s^{i_5}t^{i_6},
x^{j_1}y^{j_2}z^{j_3}]r^{j_4}s^{j_5}t^{j_6}\\
&=&
x^{j_1}y^{j_2}z^{j_3}r^{i_4+j_4}s^{i_5+j_5}t^{i_6+j_6}[r,x]^{i_4j_1}
[r,y]^{i_4j_2}[r,z]^{i_4j_3+i_6j_2}\\
&&\ \ \
{[s,x]}^{i_5j_1+i_6j_2}[s,y]^{i_5j_2}
[s,z]^{i_5j_3}[t,x]^{i_6j_1}[t,z]^{i_6j_3}.
\end{eqnarray*}
Combining these proves the lemma. \end{proof}

\begin{lemmatje} $(i_1,\cdots,i_{14})^n=$
\[ \left( \begin{array}{c}
ni_1 \\ ni_2 \\ ni_3 \\ni_4-\frac{1}{2}(n-1)ni_2i_1\\
ni_5-\frac{1}{2}(n-1)ni_3i_2\\ ni_6-\frac{1}{2}(n-1)ni_3i_1\\
ni_7+\frac{1}{2}(n-1)ni_4i_1-\frac{1}{12}(n-1)n(2n-1)i_1^2i_2
+\frac{1}{4}(n-1)ni_1i_2 \\
ni_8 +\frac{1}{2}(n-1)ni_4i_2-\frac{1}{6}(n-1)n(2n-1)i_1i_2^2
 +\frac{1}{4}(n-1)ni_1i_2-\frac{1}{4}(n-1)ni_1i_2^2\\
ni_9 +\frac{1}{2}(n-1)ni_4i_3 +\frac{1}{2}(n-1)ni_6i_2
-\frac{1}{3}(n-1)n(2n-1)i_3i_1i_2 -\frac{1}{2}(n-1)ni_1i_2i_3\\
ni_{10} +\frac{1}{2}(n-1)ni_5i_1 +\frac{1}{2}(n-1)ni_6i_2
-\frac{1}{6}(n-1)n(2n-1)i_3i_1i_2 \\
ni_{11}+\frac{1}{2}(n-1)ni_5i_2-\frac{1}{12}(n-1)n(2n-1)i_2^2i_3
+\frac{1}{4}(n-1)ni_2i_3 \\
ni_{12} +\frac{1}{2}(n-1)ni_5i_3-\frac{1}{6}(n-1)n(2n-1)i_2i_3^2
 +\frac{1}{4}(n-1)ni_3i_2-\frac{1}{4}(n-1)ni_2i_3^2\\
ni_{13}+\frac{1}{2}(n-1)ni_6i_1-\frac{1}{12}(n-1)n(2n-1)i_1^2i_3
+\frac{1}{4}(n-1)ni_1i_3 \\
ni_{14} +\frac{1}{2}(n-1)ni_6i_3-\frac{1}{6}(n-1)n(2n-1)i_1i_3^2
 +\frac{1}{4}(n-1)ni_3i_1-\frac{1}{4}(n-1)ni_1i_3^2\\

\end{array} \right) \] \end{lemmatje}

\begin{proof} By induction on $n$. \end{proof}

\begin{lemmatje} \label{berekening} Notations as above.
\[ \begin{array}{rcl}
 G^{(4),n}/G^{(4)} & \cong & \{(i_1,\ldots,i_{14})|
i_1,i_2,i_3 \in n{\msy Z},~
i_4,i_5,i_6 \in n_2{\msy Z}, \\
 & & i_7,\ldots,i_{14} \in
n_6{\msy Z},~i_9 + i_{10} \in n_2 {\msy Z} \}. \end{array} \]
(Recall that $n_i:= n / {\rm gcd} (n,i).) $\end{lemmatje}

\begin{proof}
The inclusion ``$\subset$'' follows at once from Lemma 4.2 and the following
two
observations:
$$-\frac{1}{12}(n-1)n(2n-1)i_1^2i_2+\frac{1}{4}(n-1)ni_1i_2=$$
$$\frac{1}{12}(n-1)ni_1i_2[-(2n-1)i_1+3]$$
is an element of $n_6 \msy Z$ (either $i_1$ is even or the two terms between
brackets have the same parity).
The same holds for
$$-\frac{1}{6}(n-1)n(2n-1)i_1i_2^2
 +\frac{1}{4}(n-1)ni_1i_2-\frac{1}{4}(n-1)ni_1i_2^2=$$
$$\frac{1}{12}(n-1)ni_1i_2[-(4n+1)i_2+3].$$
The other inclusion follows from a direct computation:
\[ \begin{array}{rl}
a:= & (1,0,\ldots,0)^n(-1,1,0,\ldots,0)^n(0,-1,0,\ldots,0)^n \\
 =&
 (0,0,0,\frac{1}{2}(n-1)n,0,0,-\frac{1}{6}(n-1)n(n+1),-\frac{1}{6}(n-1)n(n+1),
0,\ldots,0),\\
b:= & (-1,0,\ldots,0)^n (1,1,0\ldots,0)^n (0,-1,\ldots,0)^n\\
=& (0,0,0,-\frac{1}{2}(n-1)n,0,0,-\frac{1}{6}(n-2)(n-1)n,
\frac{1}{6}(n-1)n(n+1),0,\ldots,0),\\
c:= & (1,0\ldots,0)^n(-1,-1,0,\ldots,0)^n(0,1,0,\ldots.,0)^n \\
=&
(0,0,0,-\frac{1}{2}(n-1)n,0,0,\frac{1}{6}(n-1)n(n+1),-\frac{1}{6}(n-2)(n-1)n,
0,\ldots,0).
\end{array} \]
So we have:
$$ab=(0,0,0,0,0,0,-\frac{1}{6}(n-1)n(2n-1),0,\ldots,0),$$
$$ac=(0,0,0,0,0,0,0,-\frac{1}{6}(n-1)n(2n-1),0,\ldots,0).$$
The same computation holds for the eleventh till the fourteenth coefficient,
so it is possible to get $n_6$ as one of the coefficients
$i_7,i_8,i_{11},..,i_{14}$ and all the other coefficients zero.
Thus, looking again at $a$,
we see that it is possible to get $n_2$ as the coefficient
$i_4$ and all the other coefficients zero. By symmetry the same holds for
the coefficients $i_5,i_6$.\\
Furthermore, if we put
$$d:= (0,0,-1,0,\ldots,0)^n(0,0,1,1,0,\ldots,0)^n(0,0,0,-1,0,\ldots.,)^n,$$
$$e:= (-1,0,\ldots,0)^n(1,0,0,0,1,0,\ldots,0)^n(0,0,0,0,-1,0,\ldots.,)^n,$$
then we have:
\[ \begin{array}{l}
d=(0,0,0,0,0,0,0,0,\frac{1}{2}(n-1)n,0,\ldots,0), \\
e=(0,0,0,0,0,0,0,0,0,\frac{1}{2}(n-1)n,0,\ldots,0).
\end{array} \]
Thus, we can get $n_2$ in the
ninth and tenth coefficient, independently.
 The rest now follows from what we have proven already
and the computation of
$$(-1,0,\ldots,0)^n(0,-1,0,\ldots,0)^n(0,0,-1,0,\ldots,0)^n(1,1,1,0,\ldots,0)^n.$$
\end{proof}

\begin{cortje} \label{free comp} Notations as above. \hfill \break
1. $(G^{(2)} \cap G^{(3),n})/G^{(3)} \cong n_2G^{(2)}/G^{(3)}.$ \hfill\break
2. $(i_7,\ldots,i_{14}) \in (G^{(3)} \cap G^{(4),n})/G^{(4)} $ if and only if
$i_7,\ldots,i_{14} \in
n_6 {\msy Z},~i_9 +i_{10} \in n_2 {\msy Z}.$
\end{cortje}

\centerline{\vbox{
\hbox{M.~Pikaart and A.J.~de~Jong}
\hbox{Universiteit Utrecht,}
\hbox{Mathematisch Instituut,}
\hbox{Postbus 80.010,}
\hbox{3508 TA Utrecht,}
\hbox{The Netherlands.}
\hbox{E-mail: pikaart@math.ruu.nl \& dejong@math.ruu.nl}}}

\end{document}